\documentclass[reqno,a4paper]{amsart}
\usepackage{eucal,xspace,amssymb}
\numberwithin{equation}{section}
\newcommand{\nc}{\newcommand}
\nc{\rnc}{\renewcommand}
\nc{\DMO}{\DeclareMathOperator}
\nc{\nt}{\newtheorem}

\rnc{\thefootnote}{\alph{footnote}}

\nt{Thm}{Theorem}[section]
\nt{Lem}[Thm]{Lemma}
\nt{Prp}[Thm]{Proposition}

\theoremstyle{remark}
\nt*{Rms}{Remarks}
\nt*{acknowledgements}{\sc Acknowledgements}

\nc{\+}{^\dagger}
\nc{\0}{\mspace{-1.5mu}\circ\mspace{-0.5mu}}
\nc{\1}{\ensuremath{{\boldsymbol1}}}
\nc{\2}{{\frac{1}{2}}}
\nc{\3}{{\tfrac{1}{2}}}
\nc{\4}[1]{{\overline{#1}}}
\nc{\ABS}[1]{\lvert{#1}\rvert}
\nc{\BB}{{\mathfrak B}}
\nc{\CC}{\ensuremath{\mathbb C}\xspace}
\nc{\CK}{\ensuremath{{\mathfrak{C}(\KK)}}\xspace}
\nc{\CKK}{\ensuremath{\mathfrak{C}(\KK,\kappa)}\xspace}
\nc{\DD}{\ensuremath{\mathfrak{D}}\xspace}
\nc{\DEF}{\equiv}
\nc{\END}{\text{\rm{End}}}
\nc{\FF}{\mathfrak{F}}
\nc{\HH}{\ensuremath{\mathcal H}\xspace}
\DMO{\IND}{ind}          
\nc{\KK}{\ensuremath{\mathcal K}\xspace}                    
\nc{\NN}{\ensuremath{\mathbb N}\xspace}                  
\nc{\NORM}[1]{\lVert{#1}\rVert}
\DMO{\RAN}{ran}
\nc{\RR}{\ensuremath{\mathbb R}\xspace}
\nc{\SC}{^\sharp}
\nc{\SET}[2]{\{#1\ |\ #2\}}
\DMO{\SPAN}{span}          
\DMO{\SUPP}{supp}          
\nc{\WK}{{\WW(\KK,\kappa)}}      
\nc{\WW}{\mathfrak{W}}
\nc{\ZZ}{\ensuremath{\mathbb Z}\xspace}

\begin{document}
\title[The Charge Quantum Numbers of Quasi-free Endomorphisms]{The
Charge Quantum Numbers of Gauge Invariant Quasi-free Endomorphisms}
\author[Carsten Binnenhei]{Carsten Binnenhei\\[2ex]April 1999}
\address{{Dipartimento di Matematica, Universit\`a di Roma "La Sapienza"\\
  Piazzale Aldo Moro,~2\\I-00185 Roma}}
\email{binnenhe@mat.uniroma1.it}
\thanks{Supported by the EU TMR network ``Non-Commutative Geometry'',
contract no.\ ERB FMRX CT96--0073}
\begin{abstract}
	The representations of a group of gauge automorphisms of the
	canonical commutation or anticommutation relations which appear on
	the Hilbert spaces of isometries $H_\varrho$ implementing quasi-free
	endomorphisms $\varrho$ on Fock space are studied. Such a
	representation, which characterizes the ``charge'' of $\varrho$ in
	local quantum field theory, is determined by the Fock space
	structure of $H_\varrho$ itself: Together with a ``basic''
	representation of the group, all higher symmetric or antisymmetric
	tensor powers thereof also appear. Hence $\varrho$ is reducible
	unless it is an automorphism. It is further shown by the example of
	the massless Dirac field in two dimensions that localization and
	implementability of quasi-free endomorphisms are compatible with
	each other.
\end{abstract}
\maketitle
%%%%%%%%%%%%%%%%%%%%%%%%%%%%%%%%%%%%%%%%%%%%%%%%%%%%%%%%%%%
\section{Introduction}
In quantum field theory the structure of superselection sectors is
entirely encoded in the set of localized endomorphisms of the algebra
of local observables \cite{DHR3,H}. In the case of main physical
interest, viz.\ in four dimensional Minkowski space, the set of
(equivalence classes of) localized endomorphisms can be identified
with the representation category of a unique compact group, the
\emph{global gauge group} of the theory \cite{DR89,DR90}. This gauge
group acts on a larger field algebra containing besides the
observables charge carrying fields with normal commutation relations
which reach all superselection sectors from the vacuum
\cite{DR89a,DR90}.
%The observables can be recovered as the gauge
%invariant elements of the field algebra. 
Gauge group and field algebra
are intrinsically determined by the observable data.

The relation between localized endomorphisms and representations of
the gauge group is made concrete in the following way
\cite{DR72}. There is a functor which assigns to a localized
endomorphism $\varrho$ the \emph{Hilbert space of isometries}
$H_\varrho$ consisting of all local fields $\Psi$ which induce
$\varrho$:
$$H_\varrho\DEF\SET{\Psi}{\Psi a=\varrho(a)\Psi,\text{ for all local
observables }a}.$$ 
The action of the gauge group on the field algebra restricts to a
unitary representation $D_\varrho$ on $H_\varrho$ relative to the
inner product $\langle\Psi,\Psi'\rangle\1\DEF\Psi^*\Psi'$. This
representation of the gauge group determines the \emph{charge} of the
endomorphism $\varrho$; it is customary to refer to any label which
characterizes the representation $D_\varrho$ as the \emph{charge
quantum numbers} of $\varrho$. It will be used below that the
representation $D_\varrho$ is in a canonical way unitarily equivalent
to the representation on the Hilbert space $H_\varrho\Omega$ generated
by applying the field operators in $H_\varrho$ to the vacuum vector
$\Omega$.

Any orthonormal basis $\Psi_1,\dots,\Psi_d$ in $H_\varrho$ generates a
representation of the Cuntz algebra $\mathcal{O}_d$ and implements the
endomorphism $\varrho$ as follows:
\begin{equation}
  \label{IMP}
  \varrho(a)=\sum_{j=1}^d\Psi_ja\Psi_j^*.
\end{equation}
Using this formula, $\varrho$ can be canonically extended to an
endomorphism of the field algebra. This extension is \emph{gauge
invariant}, i.e.\ commutes with all gauge automorphisms:
\begin{Prp}
	\label{prp:GI}
	 Let $\varrho$ be an endomorphism of the field algebra which is
  implemented by a Hilbert space of isometries $H_\varrho$ as in
  \eqref{IMP}. Then $H_\varrho$ is gauge invariant if and only if
  $\varrho$ is gauge invariant.
\end{Prp}
\begin{proof}
  Assume first that $H_\varrho$ is invariant under gauge
  automorphisms $\gamma$. Since the representation $D_\varrho$ of
  the gauge group on $H_\varrho$, given by
  $D_\varrho(\gamma)\DEF\gamma|_{H_\varrho}$, 
	is unitary, the $\gamma(\Psi_j)$ also form an orthonormal basis in
  $H_\varrho$. Since the endomorphism associated with a Hilbert space
  of isometries as in \eqref{IMP} is independent of the choice of an
  orthonormal basis in $H_\varrho$, it follows that, for any field
  operator $f$,
  \begin{equation*}
		\gamma(\varrho(f))=\sum_j\gamma(\Psi_j)\gamma(f)
		\gamma(\Psi_j)^*=\varrho(\gamma(f)).
  \end{equation*}
  
  Conversely, assume that $\varrho$ is gauge invariant.  Let $\Psi\in
	H_\varrho$ and let $\gamma$ be a gauge transformation. Then one has
	for any field operator $f$
  $$\gamma(\Psi)f=\gamma\bigl(\Psi\gamma^{-1}(f)\bigr)
  =\gamma\bigl(\varrho(\gamma^{-1}(f))\Psi\bigr)=\varrho(f)\gamma(\Psi)$$
  so that $\gamma(\Psi)\in H_\varrho$.
\end{proof}

The existence of localized endomorphisms and associated Hilbert spaces
of isometries follows from first principles of local quantum field
theory. But it is by no means obvious how to obtain them explicitly in
concrete models. In previous work we have developed a general theory
of quasi-free endomorphisms of the CAR and CCR algebras which can be
implemented by Hilbert spaces of isometries on Fock space
\cite{CB1,CB3}. Among the results are implementability conditions for
endomorphisms, which generalize the well-known criteria of Shale and
Stinespring for automorphisms \cite{S,SS65}, and detailed
constructions of field operators which implement endomorphisms
according to \eqref{IMP}.

In the present paper we are interested in the possible charge quantum
numbers of such endomorphisms. The CAR resp.\ CCR algebra will play
the role of the field algebra. Therefore quasi-free endomorphisms have
to be viewed as endomorphisms of the field algebra and, by
Proposition~\ref{prp:GI}, we have to restrict attention to
endomorphisms which are gauge invariant under an appropriate group
action. We will consider quasi-free actions of arbitrary groups which
leave the Fock vacuum invariant. We show that the charge quantum
numbers are then determined by the natural Fock space structure found
in \cite{CB1,CB3} of the Hilbert spaces of isometries $H_\varrho$
implementing gauge invariant quasi-free endomorphisms $\varrho$: The
representation $D_\varrho$ is unitarily equivalent to the
representation $\Lambda_{\mathfrak{k}_\varrho}$ on the antisymmetric
Fock space over an auxiliary unitary $G$-module
$\mathfrak{k}_\varrho$, tensored with a certain character
${\det}_{\mathfrak{h}_\varrho}$:
$$D_{\varrho}\simeq{\det}_{\mathfrak{h}_\varrho}\otimes
\Lambda_{\mathfrak{k}_\varrho}\quad\text{(CAR),}$$
resp.\ to the representation $\mathfrak{S}_{\mathfrak{k}_\varrho}$ on
the symmetric Fock space over $\mathfrak{k}_\varrho$:
$$D_{\varrho}\simeq\mathfrak{S}_{\mathfrak{k}_\varrho}
\quad\text{(CCR).}$$
This is our main result, contained in Theorems~\ref{thm:CHARGE} and
\ref{ccr:thm:CHARGE}. It follows that $D_{\varrho}$ is reducible if
$\varrho$ is non-surjective. Any $\varrho$ has a quasi-free conjugate
$\varrho^c$ such that $D_{\varrho^c}$ is equivalent to the complex
conjugate of $D_{\varrho}$, provided that the single-particle space
has a particle-antiparticle symmetry.

The analysis of the representations $D_\varrho$ is completely
independent of localization properties of endomorphisms. In order to
show that localization and implementability are not in conflict, we
give in Section~\ref{sec:EXA} an explicit example of a
\emph{localized} implementable gauge invariant endomorphism $\varrho$,
with $\dim H_\varrho=2^N$, of the free massless Dirac field with
$U(N)$ gauge symmetry in two dimensions. The construction rests on the
use of ``local'' Fourier bases for the chiral components, and is in
this respect similar to the known examples \cite{JMB1,JMB2} of
localized endomorphisms in conformal field theory.

The present investigations are taken from the author's Ph.\,D. thesis
in physics \cite{CB4}, to which we refer for further results and
discussions.
%%%%%%%%%%%%%%%%%%%%%%%%%%%%%%%%%%%%%%%%%%%%%%%%%%%%%%%%%%%
\section{The Fermionic Case} \label{sec:FER}
First of all we need some formalism and some results from
\cite{CB1,CB4}.
%%%%%%%%%%%%%%%%%%%%%%%%%%%%%%%%%%%%%%%%%%%%%%%%%%%%%%%%%%%%%%%%%
\subsection{Preliminaries on the implementation of quasi-free
endomorphisms (CAR)} \label{sec:PRECAR}
Recall Araki's approach to the canonical anticommutation relations
\cite{A70,A87}: Let \KK be an infinite dimensional separable complex
Hilbert space, endowed with a complex conjugation $f\mapsto f^*$. The
{\em(selfdual) CAR algebra} \CK over \KK is the unique (simple)
$C^*$-algebra generated by \1 and the elements of \KK, subject to the
anticommutation relation
$$f^*g+gf^*=\langle f,g\rangle\1,\quad f,g\in\KK.$$
Let $P_1$ be a fixed \emph{basis projection} on \KK, i.\,e.\ an
orthogonal projection such that
$$\4{P_1}=\1-P_1.$$
Here the bar denotes the complex conjugate of an operator $A$:
$$\4{A}(f)\DEF A(f^*)^*,\quad f\in\KK.$$
Let $P_2$ be the complementary (basis) projection:
$$P_2\DEF \1-P_1.$$
The components of an operator $A$ on \KK with respect to the
decomposition
$$\KK=\KK_1\oplus\KK_2$$
given by  $P_1$ and $P_2$ will be denoted by
$$A_{mn}\DEF P_mAP_n,\quad m,n=1,2$$
and will be regarded as operators from $\KK_n$ to $\KK_m$.

To the basis projection $P_1$ there corresponds a unique (pure,
quasi-free) \emph{Fock state} $\omega_{P_1}$ which is completely
determined by the condition that
\begin{equation}
	\label{FS}
	\omega_{P_1}(f^*f)=0\quad\text{if}\quad P_1f=0.
\end{equation}
The GNS representation associated with $\omega_{P_1}$ will be denoted
by $(\FF_{P_1},\pi_{P_1},\Omega_{P_1})$.
The Hilbert space $\FF_{P_1}$ can be identified with the antisymmetric
Fock space over $\KK_1$. The elements of \KK are then represented by
sums of creation and annihilation operators:
\begin{equation}
	\label{PP1}
	\pi_{P_1}(f)= a^*(P_1f)+a(P_1f^*),\quad f\in\KK,
\end{equation}
and the cyclic vector $\Omega_{P_1}$ is the Fock vacuum vector.

Every isometry $V$ on \KK which commutes with complex conjugation
extends to a unique \emph{quasi-free endomorphism} $\varrho_V$ of \CK:
$$\varrho_V(f)=V(f),\quad f\in\KK.$$
As shown in \cite{CB1}, an endomorphism $\varrho_V$ can be implemented
by a Hilbert space of isometries $H_{\varrho_V}$ (cf.~\eqref{IMP}) in
the Fock state $\omega_{P_1}$ if and only if  
\begin{equation}
	\label{IMPCOND}
	[P_1,V]\text{ is Hilbert--Schmidt (HS)}.
\end{equation}
These isometries form a semigroup
$$\END_{P_1}(\KK)\DEF\SET{V\in\BB(\KK)}{V^*V=\1,\ \4{V}=V,\ [P_1,V]
\text{ is HS}}$$
isomorphic to the semigroup of all implementable quasi-free
endomorphisms of \CK. The \emph{statistics dimension} $d_{\varrho_V}$
of $\varrho_V$ is given by
$$d_{\varrho_V}\DEF\dim H_{\varrho_V}=2^{\2\IND V},$$
where $\IND V$ is, up to the sign, the Fredholm index of $V$:
$$\IND V=\dim\ker V^*\in\{0,2,4,\dots,\infty\},\quad
V\in\END_{P_1}(\KK).$$

The grading automorphism of \CK is equal to the quasi-free
automorphism $\varrho_{-\1}$. Let $\Gamma(-\1)$ be the (self-adjoint,
unitary) second quantization of $\varrho_{-\1}$, given by
$$\Gamma(-\1)\pi_{P_1}(a)\Omega_{P_1}=\pi_{P_1}(\varrho_{-\1}(a))
\Omega_{P_1},\quad a\in\CK,$$
and let $\theta(-\1)$ be the unitary operator
\begin{equation}
	\label{THETA}
	\theta(-\1)\DEF\tfrac{1}{\sqrt{2}}\bigl(\1-i\Gamma(-\1)\bigr).
\end{equation}
Then the \emph{twisted Fock representation} $\psi_{P_1}$ induced by
$P_{1}$ is defined by
\begin{equation}
	\label{TWIST}
	\psi_{P_1}(a)\DEF\theta(-\1)\pi_{P_1}(a)\theta(-\1)^*.
\end{equation}
It can be used to describe the commutants of ``local'' subalgebras:
If $\HH\subset\KK$ is a subspace invariant under complex conjugation,
and $\mathfrak{C}(\HH)$ the $C^*$-subalgebra of \CK generated by
$\HH$, then
$$\pi_{P_1}(\mathfrak{C}(\HH))'=\psi_{P_1}(\mathfrak{C}(\HH^\bot))''
\qquad\text{(\emph{twisted duality})}.$$

Let $V\in\END_{P_1}(\KK)$ be given. As mentioned in the introduction,
it suffices for our purposes to consider group actions on the Hilbert
space $H_{\varrho_V}\Omega_{P_1}$. An orthonormal basis in this space
can be obtained as follows. Define a finite dimensional subspace
$\mathfrak{h}\subset\KK_1$ by
\begin{equation}
	\label{H}
	\mathfrak{h}\DEF V_{12}(\ker V_{22}),
\end{equation}
and an antisymmetric\footnote{An operator $A$ on \KK is
\emph{antisymmetric} if $\4{A^*}=-A$.} Hilbert--Schmidt operator $T$
from $\KK_1$ to $\KK_2$ by
\begin{equation}
	\label{T}
	T\DEF V_{21}{V_{11}}^{-1}-{V_{22}}^{-1*}{V_{12}}^*[\ker {V_{11}}^*].
\end{equation}
Here the bounded operator ${V_{11}}^{-1}$ is defined to be zero on
$\ker{V_{11}}^*$, and $[\HH]$ denotes the orthogonal projection onto a
closed subspace $\HH\subset\KK$. Then one has $T(\mathfrak{h})=0$;
such pairs $(\mathfrak{h},T)$ parameterize the class of all Fock
states which are unitarily equivalent to the given Fock state
$\omega_{P_1}$ \cite {CB4}. The basis projection $P$ corresponding to
the pair $(\mathfrak{h},T)$ is explicitly given by 
\begin{equation}
	\label{P}
	P\DEF(P_1+T)(P_1+T^*T)^{-1}(P_1+T^*)-[\mathfrak{h}]+
	[\mathfrak{h}^*],
\end{equation}
and $\mathfrak{h}$ and $T$ can be recovered from $P$ by
\begin{align}
	\mathfrak{h} &=\ker P_{11},\label{HDEF}\\
	T &=P_{21}{P_{11}}^{-1}\label{TDEF}
\end{align}
(${P_{11}}^{-1}$ is defined in a similar way as ${V_{11}}^{-1}$ above).
The Fock state $\omega_P$ associated with $P$ is an extension of the
partial Fock state $\omega_{P_1}\0\varrho_V^{-1}|_{\RAN\varrho_V}$,
and is induced by the cyclic vector
$$\Omega_P=\bigl({\det}(P_1+T^*T)\bigr)^{-1/4}\psi_{P_1}(e_1\dotsm
e_L)\exp(\3\4{T}a^*a^*)\Omega_{P_1}\in\FF_{P_1}$$ 
where the determinant has to be computed on $\KK_1$,
$\{e_1,\dots,e_L\}$ is an orthonormal basis in $\mathfrak{h}$, and the
exponential term has a well-defined meaning as a strongly convergent
series on a dense domain containing $\Omega_{P_1}$. The vector
$\Omega_P$ belongs to $H_{\varrho_V}\Omega_{P_1}$; in fact, the latter
Hilbert space consists precisely of the vectors in $\FF_{P_1}$ which
induce extensions of the partial Fock state
$\omega_{P_1}\0\varrho_V^{-1}|_{\RAN\varrho_V}$.

A complete orthonormal basis in $H_{\varrho_V}\Omega_{P_1}$ can be
obtained by applying suitable partial isometries from the commutant of
$\RAN\varrho_V$ to $\Omega_P$. The basis projection $P$ leaves $\ker
V^*$ invariant. Let
\begin{equation}
	\label{K}
	\mathfrak{k}\DEF P(\ker V^*),\quad\text{with }
	\dim\mathfrak{k}=\3\IND V,
\end{equation}
and let $\{g_j\}$ be an orthonormal basis in $\mathfrak{k}$. For any
multi-index $\alpha=(\alpha_1,\dots,\alpha_l),\
1\leq\alpha_1<\dots<\alpha_l\leq\2\IND V$ (resp.\ $\alpha=0$ if
$l=0$), set 
\begin{equation}
	\label{OA}
	\Omega_\alpha\DEF\psi_{P_1}(g_{\alpha_1}\dotsm
	g_{\alpha_l})\Omega_P.
\end{equation}
Then the $\Omega_\alpha$ constitute an orthonormal basis in
$H_{\varrho_V}\Omega_{P_1}$, and they determine an orthonormal basis
$\{\Psi_\alpha\}$ in $H_{\varrho_V}$ via
\begin{equation}
	\label{HRV}
	\Psi_\alpha\pi_{P_1}(a)\Omega_{P_1}=
	\pi_{P_1}(\varrho_V(a))\Omega_\alpha.
\end{equation}
Since $\psi_{P_1}$ is a representation of the canonical
anticommutation relations and since $\Omega_P$ is annihilated by the
operators $\psi_{P_1}(g_j)^*$, the spaces $H_{\varrho_V}$ and
$H_{\varrho_V}\Omega_{P_1}$ can both be identified with the
antisymmetric Fock space over $\mathfrak{k}$ \cite{CB1,CB4}. 
%%%%%%%%%%%%%%%%%%%%%%%%%%%%%%%%%%%%%%%%%%%%%%%%%%%%%%%%%%%%%%%%%
\subsection{Gauge invariant quasi-free endomorphisms (CAR)}
\label{sec:GICAR}
Let $G\subset\END_{P_1}(\KK)$ be a group\footnote{The results in this
section hold also for non-compact groups (relative to the strong
topology). However, the close relationship between representations of
$G$ and superselection sectors is then lost. We also do not require
that $-\1\in G$, which would be necessary if the $G$-invariant
elements of \CK were to be interpreted as physical observables.}
consisting of unitary operators which commute with $P_1$. The usual
second quantization of $U\in G$ (or, more precisely, of $U_{11}$) will
be denoted by $\Gamma(U)$; the map $U\mapsto\Gamma(U)$ is strongly
continuous. The corresponding gauge automorphisms, which leave
$\omega_{P_1}$ invariant, will be denoted by $\gamma_U$.

We are interested in the representation $D_{\varrho_V}$ of $G$ on the
Hilbert space $H_{\varrho_V}$ which implements a gauge invariant
quasi-free endomorphism $\varrho_V$. \emph{Gauge invariant}
implementable quasi-free endomorphisms are given by the elements of
the semigroup
$$\END_{P_1}(\KK)^G\DEF\SET{V\in\END_{P_1}(\KK)}{[V,G]=0}.$$ 
To determine $D_{\varrho_V}$ up to unitary equivalence, it suffices to
calculate the transformed vectors $\Gamma(U)\Omega_\alpha$, $U\in G$,
where the $\Omega_\alpha$ are the basis elements in
$H_{\varrho_V}\Omega_{P_1}$ defined in \eqref{OA}.  The basic
observation is that the objects entering the construction of the
$\Omega_\alpha$, namely the spaces $\mathfrak{h}$ and $\mathfrak{k}$
and the operator $T$, are all gauge invariant.
\begin{Lem}
  \label{lem:GP} 
	Let $U\in G$. Then $\Gamma(U)$ implements the gauge automorphism
	$\gamma_U$ in the twisted Fock representation $\psi_{P_1}$:
  $$\Gamma(U)\psi_{P_1}(a)=\psi_{P_1}(\gamma_U(a))\Gamma(U),
	\quad a\in\CK.$$
\end{Lem}
\begin{proof}  
	Since $U$ commutes with $P_1$, the implementer $\Gamma(U)$ is even:
  $$[\Gamma(U),\Gamma(-\1)]=0.$$
  This implies that $[\Gamma(U),\theta(-\1)]=0$ (see \eqref{THETA}) so
	that, by \eqref{TWIST},
	\begin{equation*}
		\Gamma(U)\psi_{P_1}(a)\Gamma(U)^*=\theta(-\1)\Gamma(U)
		\pi_{P_1}(a)\Gamma(U)^*\theta(-\1)^*=\psi_{P_1}(\gamma_U(a)).
	\end{equation*}
\end{proof}
\begin{Lem}
  \label{lem:EXP}
  Let $V\in\END_{P_1}(\KK)^G$. Then $\exp(\3\4{T}a^*a^*)\Omega_{P_1}$,
  with $T$ defined by \eqref{T}, is invariant under all gauge
  transformations $\Gamma(U),\ U\in G$. 
\end{Lem}
\begin{proof}
  If $S$ is an antisymmetric operator from $\KK_1$ to $\KK_2$ of
  finite rank, then one readily verifies that
  $$\Gamma(U)(\3\4{S}a^*a^*)\Gamma(U)^*=\3(U\4{S}U^*)a^*a^*$$ 
	(here the expressions of the form $\4{S}a^*a^*$ are defined by
	expanding $S=\sum f_j^*\langle g_j,.\rangle$ with $f_j,g_j\in\KK_1$,
	and by setting $\4{S}a^*a^*=\sum a^*(f_j)a^*(g_j)$).  Approximating
	$T$ by such finite rank operators in the Hilbert--Schmidt norm
	(cf.~\cite{CR}), one finds that  
  $$\Gamma(U)(\3\4{T}a^*a^*)^n\Omega_{P_1}=
  \bigl(\3(U\4{T}U^*)a^*a^*\bigr)^n\Omega_{P_1},\quad n\in\NN,$$
  because $\Gamma(U)\Omega_{P_1}=\Omega_{P_1}$. It follows that
  \begin{equation*}
    \begin{split}
      \Gamma(U)\exp(\3\4{T}a^*a^*)\Omega_{P_1}&=\sum_{n=0}^\infty
        \frac{1}{n!}\Gamma(U)(\3\4{T}a^*a^*)^n\Omega_{P_1}\\
      &=\exp\bigl(\3(U\4{T}U^*)a^*a^*\bigr)\Omega_{P_1}.
    \end{split}
  \end{equation*}
  Since $U$ commutes with $P_1,\ P_2$ and $V$, it also commutes with
  all components of $V$ and $V^*$, including the operators
  ${V_{11}}^{-1},\ [\ker{V_{11}}^*]$ etc. The operator $T$ is by
  \eqref{T} a bounded function of these components, so that
  \begin{equation}
    \label{GT}
    [U,T]=0,\quad U\in G.
  \end{equation}
  Hence we get
  $\Gamma(U)\exp(\3\4{T}a^*a^*)\Omega_{P_1}=
  \exp(\3\4{T}a^*a^*)\Omega_{P_1}$ as claimed.
\end{proof}
Setting $\delta\DEF\bigl({\det}(P_1+T^*T)\bigr)^{-1/4}$, we thus arrive
at the following formula:
\begin{equation}
  \label{GPSIOM}
  \Gamma(U)\Omega_\alpha=\delta\cdot\psi_{P_1}\bigl(U(g_{\alpha_1})
	\dotsm U(g_{\alpha_l})U(e_1)\dotsm U(e_L)\bigr)\exp(\3\4{T}a^*a^*)
	\Omega_{P_1}.
\end{equation}
%which enables us to describe the representation $D_{\varrho_V}$ of $G$
%on $H_{\varrho_V}$.
\begin{Thm}
  \label{thm:CHARGE} 
	Let $P_1$ be a basis projection of \KK, let $G$ be a group of
  unitary operators on \KK commuting with $P_1$ and with complex
  conjugation, and let $V\in\END_{P_1}(\KK)^G$. Then the finite
  dimensional subspace $\mathfrak{h}\subset\KK_1$ and the $(\2\IND
  V)$ dimensional subspace $\mathfrak{k}\subset\KK$ associated with
  $V$ by \eqref{H} and \eqref{K} are both invariant under $G$. Let
  $\Lambda_{\mathfrak{k}}$ be the unitary representation of $G$ on the
  antisymmetric Fock space over $\mathfrak{k}$ that is obtained by
  taking antisymmetric tensor powers of the representation on
  $\mathfrak{k}$. Then the unitary representation $D_{\varrho_V}$ of
  $G$ on the Hilbert space of isometries $H_{\varrho_V}$ which
  implements $\varrho_V$ in the Fock state $\omega_{P_1}$ is unitarily
  equivalent to $\Lambda_{\mathfrak{k}}$, tensored with the
  one dimensional representation
  ${\det}_{\mathfrak{h}}(U)\DEF\det(U|_{\mathfrak{h}})$\textup{:}
  \begin{equation} 
		\label{REPG}
  	D_{\varrho_V}\simeq{\det}_{\mathfrak{h}}\otimes
  	\Lambda_{\mathfrak{k}}.
  \end{equation}
\end{Thm}
\begin{proof}
  The subspace $\mathfrak{h}=V_{12}(\ker V_{22})$ is invariant under
  $G$ because $G$ commutes with the components of $V$ (cf.\ the proof
  of Lemma~\ref{lem:EXP}). Since $\{e_1,\dots,e_L\}$ is an orthonormal
  basis in $\mathfrak{h}$, it follows from the canonical
  anticommutation relations that
  $$U(e_1)\dotsm U(e_{L})=\det\bigl(U|_{\mathfrak{h}}\bigr)\cdot
  e_1\dotsm e_{L},\qquad U\in G.$$
  
  Similarly, $G$ commutes with the basis projection $P$
  (cf.\ \eqref{P} and \eqref{GT}) and leaves $\ker V^*$ invariant, so
  that $\mathfrak{k}=P(\ker V^*)$ is also left invariant. It then
  follows from the canonical anticommutation relations that
  $g_{\alpha_1}\dotsm g_{\alpha_l}$ transforms like the $l$-fold
  antisymmetric tensor product of $g_{\alpha_1},\dots,g_{\alpha_l}$
  under $G$.
  
  Thus we see from \eqref{GPSIOM} that the representation of $G$ on
  $H_{\varrho_V}\Omega_{P_1}$ is unitarily equivalent to
  ${\det}_{\mathfrak{h}}\otimes\Lambda_{\mathfrak{k}}$, and the same
  holds true for the representation $D_{\varrho_V}$.
\end{proof}
\begin{Rms}
  (i) Theorem~\ref{thm:CHARGE} shows that non-surjective quasi-free
	endomorphisms $\varrho_V$ are always \emph{reducible} in the sense
	that the representation $D_{\varrho_V}$ (or, if $G$ is compact, the
	representation induced by $\varrho_V$ of the pointwise gauge
	invariant ``observable'' algebra $\CK^G$ on the subspace of
	$\Gamma(G)$-invariant vectors in $\FF_{P_1}$) is reducible. In fact,
	each ``$n$-particle'' subspace of $H_{\varrho_V}$, i.e.\ the closed
	linear span of all $\Psi_\alpha$ with $\alpha$ of length $n$, is
	invariant under $G$, and may decompose further.  Let
	$D_{\varrho_V}^{(n)}$ be the restriction of $D_{\varrho_V}$ to this
	subspace. Closest to irreducibility is the case that at least
	$D_{\varrho_V}^{(1)}$ is irreducible. In typical situations, the
	remaining representations $D_{\varrho_V}^{(n)}$ will then also be
	irreducible. This happens for instance if $G\cong U(N)$ or $G\cong
	SU(N)$, and $\mathfrak{k}$ carries the defining representation of
	$G$. In the $U(N)$ case, the $D_{\varrho_V}^{(n)}$ are not only
	irreducible, but also mutually inequivalent. In the $SU(N)$ case,
	the representations
	$D_{\varrho_V}^{(0)},\dots,D_{\varrho_V}^{(N-1)}$ are mutually
	inequivalent, but $D_{\varrho_V}^{(N)}$ is equivalent to
	$D_{\varrho_V}^{(0)}$. In general, it can nevertheless happen that
	$D_{\varrho_V}^{(1)}$ is irreducible but some $D_{\varrho_V}^{(n)}$
	are not, as is the case if $G\cong SO(N)$ ($N>2$ even) and
	$\mathfrak{k}$ carries the defining representation of $G$ (cf.\
	\cite{Wey,Boe}). If already $D_{\varrho_V}^{(1)}$ is reducible, then
	one has an additional Clebsch--Gordan type splitting.

	(ii) Theorem~\ref{thm:CHARGE} characterizes the representation
  $D_{\varrho_V}$ associated with a fixed gauge invariant endomorphism
  $\varrho_V$ in terms of the representations on $\mathfrak{h}$ and
  $\mathfrak{k}$. The question which representations of $G$ can be
  realized on the spaces $\mathfrak{h}$ and $\mathfrak{k}$ by letting
  $V$ vary through $\END_{P_1}(\KK)^G$ is studied in \cite{CB4}. In
  typical field theoretic situations, where the single-particle space
  $\KK_1$ carries an irreducible representation of a compact group
  together with its complex conjugate, both with infinite
  multiplicity, as e.\,g.\ in Section~\ref{sec:EXA} below, one finds
  that any irreducible representation of $G$ realized on Fock space
  $\FF_{P_1}$ is equivalent to a subrepresentation of some
  $D_{\varrho_V}$.

	(iii) A special case worth mentioning is the case $G\cong U(1)$ and
	$\IND V=0$, i.e.\ the case of the \emph{restricted unitary group.}
	It is well-known from the work on the external field problem (see
	e.\,g.\ \cite{CHOB}) that the charge of elements of the restricted
	unitary group is given by a certain Fredholm index $\IND V_{++}$
	(which has nothing to do with the index of $V$, but refers to a
	finer decomposition $V_{11}=V_{++}\oplus V_{--}$). This fact can be
	easily derived from our general result: The factor
	$\Lambda_{\mathfrak{k}}$ in \eqref{REPG} becomes trivial, whereas
	$${\det}_{\mathfrak{h}}(U_\lambda)=\exp(i\lambda\IND V_{++})$$
	if $U_\lambda\in G$ corresponds to $e^{i\lambda}\in U(1)$ \cite{CB4}.

	Similarly, in the case $G=\{\pm\1\}\cong\ZZ_2$ and $\IND V=0$, the
	factor $\Lambda_{\mathfrak{k}}$ is trivial, but
	$${\det}_{\mathfrak{h}}(-\1)=(-1)^{\dim\mathfrak{h}}=
	(-1)^{\dim\ker V_{11}}$$
  yields the $\ZZ_2$-index of Araki and Evans \cite{AE,A87,EK}.

	(iv) If the single-particle space $\KK_1$ decomposes into the direct
	sum of two anti\-unitarily equivalent $G$-modules
	(``particle-antiparticle symmetry''), then there exists an
	involutive automorphism $V\mapsto V^c$ of $\END_{P_1}(\KK)^G$ such
	that the spaces $\mathfrak{h}^c$ and $\mathfrak{k}^c$ corresponding
	by \eqref{H} and \eqref{K} to $V^c$ are, as $G$-modules, antiunitarily
	equivalent to the spaces $\mathfrak{h}$ and $\mathfrak{k}$
	corresponding to $V$. That is, the representation
	$D_{\varrho_{V^c}}$ is unitarily equivalent to the complex conjugate
	of $D_{\varrho_{V}}$ (\emph{charge conjugation}).
\end{Rms}
%%%%%%%%%%%%%%%%%%%%%%%%%%%%%%%%%%%%%%%%%%%%%%%%%%%%%%%%%%%%%%%%%
\subsection{An example: Localized endomorphisms of the chiral Dirac
	field} \label{sec:EXA}
In Section~\ref{sec:GICAR} we have analyzed the charge quantum numbers
of gauge invariant implementable quasi-free endomorphisms in complete
generality. In particular, and in sharp contrast to the field
theoretic situation, it was not necessary to assume any localization
properties of endomorphisms. If one could find, in a specific model, a
localized implementable quasi-free endomorphism, then our methods
would apply and could be used to determine its charge and to construct
the corresponding local fields. It is however not clear from the
outset whether localization and implementability are compatible with
each other\footnote{Known results concerning this question are
restricted to the case of automorphisms. Building on the work of Carey
and Ruijsenaars \cite{CR} and others, we constructed in \cite{CB0} a
family of (implementable and transportable) localized automorphisms,
carrying arbitrary $U(1)$-charges, of the free Dirac field in two
spacetime dimensions with arbitrary mass. The operators
$V\in\END_{P_1}(\KK)^{U(1)}$ belonging to these automorphisms are
given by two $U(1)$-valued functions which are equal to $1$ at
spacelike infinity, and the charge $\IND V_{++}$ (cf.\ Remark~(iii) in
Section~\ref{sec:GICAR}) of $\varrho_V$ is equal to the difference of
the winding numbers of these functions. However, unlike in two
dimensions, there seem to be no known examples of implementable
charge-carrying automorphisms in the case $G\cong U(1)$ in four
spacetime dimensions.}. To show that this is in fact the case, we will
present below an explicit example of a non-surjective implementable
localized quasi-free endomorphism of the free massless Dirac field in
two spacetime dimensions.

Let us first introduce the free Dirac field with global $U(N)$
symmetry.  Let
\begin{equation}
	\label{HL2}
	\HH\DEF L^2(\RR^{2n-1},\CC^{2^n})
\end{equation}
be the single-particle space of the time-zero Dirac field in $2n$
spacetime dimensions. Let $H=-i\vec\alpha\vec\nabla+\beta m$ be the
free Dirac Hamiltonian, with spectral projections $E_\pm$
corresponding to the positive resp.\ negative part of the spectrum of
$H$. Tensored with $\1_N$, these operators act on the space
\begin{equation}
	\label{CN}
	\HH'\DEF\HH\otimes\CC^N.
\end{equation}
The gauge group $U(N)$ also acts naturally on $\HH'$. In the selfdual
CAR formalism, one sets
\begin{equation}
  \label{H'}
  \KK\DEF\HH'\oplus{\HH'}^*,
\end{equation}
where ${\HH'}^*$ is the Hilbert space conjugate to $\HH'$. There is a
natural conjugation $f\mapsto f^*$ on \KK which is inherited from the
antiunitary identification map $\HH'\to{\HH'}^*$. The basis projection
$P_1$ corresponding to the vacuum representation of the field is given
by
$$P_1\DEF E_+'\oplus\4{E_-'}$$
with $E_\pm'=E_\pm\otimes\1_N$. Gauge transformations act like
$U=(\1_\HH\otimes u)\oplus(\1_{\HH^*}\otimes\4{u}),\ u\in U(N)$, on
\KK. They commute with $P_1$. The field operators $\varphi_t$ at time
$t$ are given by
$$\varphi_t(f)\DEF\pi_{P_1}(e^{itH'}f)=
a(E_+'e^{itH'}f)^*+a(\4{E_-'}e^{-it\4{H'}}f^*)$$
with $H'\DEF H\otimes\1_N,\ f\in\HH'$. They are solutions of the
Dirac--Schr\"odinger equation
$$-i\frac{d}{dt}\varphi_t(f)=\varphi_t(H'f),\quad f\in D(H').$$
If $O$ is a double cone in Minkowski space with base
$B\subset\RR^{2n-1}$ at time $t$, then the local field algebra
associated with $O$ is the von Neumann algebra generated by all
$\varphi_t(f)$ with $\SUPP f\subset B$. The local observable algebras
are the fixed point subalgebras of the local field algebras under the
gauge action. A whole net of local algebras is generated from these
special ones by applying Lorentz transformations.

Gauge invariant implementable localized endomorphisms of the
$N$-component Dirac field can be characterized as follows. A
quasi-free endomorphism $\varrho_V$ is gauge invariant if and only if
$V$ has the form
\begin{equation}
  \label{VFORM}
  V=(v\otimes\1_N)\oplus(\4{v}\otimes\4{\1_N})
\end{equation}
with respect to the decomposition \eqref{H'}, where $v$ is an isometry
of \HH. This follows from the fact that the defining representation of
$U(N)$ and its complex conjugate are disjoint, so that the commutant
of $G$ on \KK is given by
$$G'=\bigl(\BB(\HH)\otimes\1_N\bigr)\oplus\bigl(\BB(\HH^*)\otimes
\4{\1_N}\bigr).$$
For $V$ of the form \eqref{VFORM} one has
$$[P_1,V]=\Bigl([E_+,v]\otimes\1_N\Bigr)\oplus\Bigl(\4{[E_-,v]}
\otimes\4{\1_N}\Bigr)$$
so that the implementability condition \eqref{IMPCOND} holds if and
only if
\begin{equation}
	\label{Ev}
	[E_+,v]\text{ and }[E_-,v]\text{ are Hilbert--Schmidt}.
\end{equation}
Therefore $\END_{P_1}(\KK)^{U(N)}$ is isomorphic to the semigroup of
all isometries $v$ of \HH which fulfill \eqref{Ev}.

An endomorphism of the algebra of all local observables is
\emph{localized} in a bounded region $O$ in Minkowski space if it acts
like the identity on observables which are localized in bounded
regions contained in the spacelike complement $O'$ of $O$
\cite{H}. Localized elements of $\END_{P_1}(\KK)^{U(N)}$ (at time
zero) can be characterized as follows.
\begin{Prp}
  \label{prop:LOC}
  Let $O$ be a double cone with base $B\subset\RR^{2n-1}$ at time
  zero.  Let $V$ be an element of $\END_{P_1}(\KK)^{U(N)}$, and let
  $v$ be the isometry of \HH associated with $V$ by
  \eqref{VFORM}. Then $\varrho_V$ is localized\footnote{More
  precisely, the normal extension of $\varrho_V$ in the representation
  $\pi_{P_1}$ is localized in $O$.}  in $O$ if and only if there
  exists, for each connected component $\Delta$ of
  $\RR^{2n-1}\setminus B$, a phase factor $\tau_\Delta\in{U(1)}$ such
  that \begin{equation} \label{TAUDELTA} v(f)=\tau_\Delta
  f\quad\text{for all }f\in\HH\text { with } \SUPP f\subset\Delta.
  \end{equation}
\end{Prp}
\begin{proof}
  Assume that $\varrho_V$ is localized in $O$. Let $b_1,\dots,b_N$ be
  the standard basis in $\CC^N$, let $\Delta$ be a component of the
  complement of $B$, and let $f,g\in\HH$ with $\SUPP
  f,g\subset\Delta$. Then
  $$a(f,g)\DEF\sum_{j=1}^N(f\otimes b_j)(g\otimes b_j)^*$$
  is gauge invariant, and $\pi_{P_1}(a(f,g))$ is an observable
  localized in $O'$. Since $\varrho_V$ is localized in $O$, one has
  $a(f,g)=\varrho_V(a(f,g))=\sum_j(v(f)\otimes b_j)(v(g)\otimes b_j)^*$.
  Since the $b_j$ are linearly independent, it follows that
  \begin{equation}
    \label{FG}
    (f\otimes b_j)(g\otimes b_j)^*=(v(f)\otimes b_j)(v(g)\otimes b_j)^*,
    \qquad j=1,\dots,N.
  \end{equation}
  Now let $P'$ be the (basis) projection onto $\HH'\subset\KK$, and let
	$\omega_{P'}$ be the corresponding Fock state. One has 
  $$\omega_{P'}\bigl((f\otimes b_j)^*(f\otimes b_j)(f\otimes b_j)^*
  (f\otimes b_j)\bigr)=\NORM{f}^4,$$
  and, since $(v(f)\otimes b_j)^*$
  belongs to the annihilator ideal of $\omega_{P'}$,
	\begin{equation*}
  	\omega_{P'}\bigl((f\otimes b_j)^*(v(f)\otimes b_j)(v(f)
		\otimes b_j)^*(f\otimes b_j)\bigr)=\ABS{\langle v(f),f\rangle}^2.
	\end{equation*}
  Therefore one gets from \eqref{FG}, in the special case $f=g$, that
  $\NORM{f}^2=\ABS{\langle v(f),f\rangle}$. It follows that there
  exists $\tau_f\in{U(1)}$ such that $v(f)=\tau_f f$. By the same
  argument, $v(g)=\tau_g g$ for some $\tau_g\in{U(1)}$. Then \eqref{FG}
  yields that $\tau_f=\tau_g$. Therefore these phase factors depend
  only on $\Delta$ and not on the functions.
  
  Conversely, assume that \eqref{TAUDELTA} holds. Then $\varrho_V$
  acts on fields localized in bounded regions in $O'$ like a gauge
  transformation, and therefore like the identity on observables
  localized in $O'$.  It follows that $\varrho_V$ is localized in $O$.
\end{proof}
Of course, $\RR^{2n-1}\setminus B$ is connected if $n>1$, but it has
two connected components if $n=1$. This is the basic reason for the
possible occurrence of braid group statistics and soliton sectors in
two dimensional Minkowski space.

Next let us demonstrate that at least the free massless Dirac field in
two spacetime dimensions possesses non-surjective implementable
localized quasi-free endomorphisms. It suffices to consider one chiral
component of the field. Thus consider the Hilbert space $\HH=L^2(\RR)$
with Dirac Hamiltonian $-i\tfrac{d}{dx}$. It is convenient to
transform to the Hilbert space $L^2({S^1})$ via the Cayley transform
$\vartheta$
$$\vartheta: \RR\cup\{\infty\}\to{S^1},\qquad x\mapsto -e^{2i\arctan
  x}=\frac{x-i}{x+i}$$
(cf.\ \cite{CR}). $\vartheta$ induces a unitary transformation
$\tilde\vartheta$ 
$$\tilde\vartheta:L^2({S^1})\to L^2(\RR),\qquad(\tilde\vartheta f)(x)=
\pi^{-\2}\frac{f(\vartheta(x))}{x+i}.$$
The important point is that the spectral projections $E_\pm$ of
$-i\tfrac{d}{dx}$ are transformed into the Hardy space projections:
Set $\tilde E_\pm\DEF\tilde\vartheta^{-1}E_\pm\tilde\vartheta$, then
$$\tilde E_+=\sum_{n\geq0}e_n\langle e_n,.\,\rangle,\quad 
\tilde E_-=\sum_{n<0}e_n\langle e_n,.\,\rangle,\quad e_n(z)\DEF z^n\
(z\in{S^1},n\in\ZZ).$$
We want to construct an isometry $v$ of $L^2({S^1})$ with
$[\tilde E_\pm,v]$ Hilbert--Schmidt (implementability), with $\IND v=1$
(close to irreducibility, cf.\ Rem.~(i) in Section~\ref{sec:GICAR}),
and such that $v(f)=f$ for all $f\in L^2({S^1})$ with $\SUPP
f\subset{S^1}\setminus I$, where $I\subset{S^1}$ is a fixed interval
(localization). As localization region we shall choose the interval
$$I\DEF\SET{e^{i\lambda}}{\tfrac{\pi}{2}\leq\lambda\leq\tfrac{3\pi}{2}}$$
which corresponds, by the inverse Cayley transform, to the interval
$\vartheta^{-1}(I)=[-1,1]$ in \RR. We need the following orthonormal
basis $(f_m)_{m\in\ZZ}$ in $L^2(I)\subset L^2({S^1})$
$$f_m(z)\DEF\sqrt{2}(-1)^mz^{2m}\chi_I(z),\quad z\in{S^1},$$
where $\chi_I$ is the characteristic function of $I$. We now define
the isometry $v$ by
\begin{equation}
	\label{v}
	v\DEF\1+\sum_{m\geq0}(f_{m+1}-f_m)\langle f_m,.\,\rangle.
\end{equation}
Note that $v$ acts like the identity on functions with support in
${S^1}\setminus I$, that $v(f_m)=f_m$ if $m<0$, and that $v$ acts like the
unilateral shift on the remaining $f_m$: $v(f_m)=f_{m+1}$ if $m\geq0$.
\begin{Lem}
  \label{lem:HS}
  The commutators $[\tilde E_+,v]$ and $[\tilde E_-,v]$ are Hilbert--Schmidt.
\end{Lem}
\begin{proof}
	The rather lengthy estimates of the Hilbert--Schmidt norms of these
	commutators, which are essentially due to P.~Grinevich, can be found
	in \cite{CB4}.
\end{proof}
Thus the operator $V\in\END_{P_1}(\KK)^{U(N)}$ induced by $v$ via
\eqref{VFORM} yields a localized endomorphism $\varrho_V$ of the
chiral Dirac field. Since by construction
$$\3\IND V=N,$$
it is clear that the space $\mathfrak{k}$ associated with $V$ by
\eqref{K} carries either the defining representation of $U(N)$ or its
complex conjugate.  By Remark~(i) in Section~\ref{sec:GICAR}, the
irreducible constituents of $\varrho_V$ correspond to the irreducible,
mutually inequivalent representations $D_{\varrho_V}^{(n)},\
n=0,\dots,N$. Note that the same isometry $v$ gives rise to localized
gauge invariant implementable endomorphisms for \emph{arbitrary}
symmetry groups $G$, by replacing the defining representation of
$U(N)$ in Eq.~\eqref{CN} with a suitable finite dimensional
representation of $G$.
%%%%%%%%%%%%%%%%%%%%%%%%%%%%%%%%%%%%%%%%%%%%%%%%%%%%%%%%%%%%%%%%%
\section{The Bosonic Case}
We need some preparations from \cite{CB3}. The exposition will closely
follow the lines of the Fermionic case considered in
Section~\ref{sec:FER}. 
%%%%%%%%%%%%%%%%%%%%%%%%%%%%%%%%%%%%%%%%%%%%%%%%%%%%%%%%%%%%%%%%%
\subsection{Preliminaries on the implementation of quasi-free
endomorphisms (CCR)} \label{sec:PRECCR}
We start with a Fock representation of the canonical commutation
relations. Thus we may assume as above that \KK is an infinite
dimensional separable complex Hilbert space with a complex conjugation
$f\mapsto f^*$, and that $P_1$ is a fixed basis projection. Let $P_2$
be the complementary projection. Define a self-adjoint unitary
operator
$$C\DEF P_1-P_2$$
so that $\4{C}=-C$, and a nondegenerate hermitian sesquilinear form
$$\kappa(f,g)\DEF\langle{f,Cg}\rangle$$
so that
$$\kappa(f^*,g^*)=-\kappa(g,f),\quad f,g\in\KK.$$
It must be emphasized that the basic form on \KK which
determines the canonical commutation relations is $\kappa$ and not the
Hilbert space inner product. In fact, the latter depends on the choice
of the Fock state, i.\,e.\ on the choice of $P_1$.

The \emph{(selfdual) CCR algebra} \CKK over $(\KK,\kappa)$ is the
simple *-algebra which is generated by \1 and elements $f\in\KK$,
subject to the commutation relation \cite{AS71,A71}
\begin{equation*}
  f^*g-gf^*=\kappa(f,g)\1,\quad f,g\in\KK.
\end{equation*}
The \emph{Weyl algebra} $\WK$ over $(\KK,\kappa)$ is the simple
$C^*$-algebra generated by unitary operators $w(f),\ f\in\KK$ with
$f=f^*$, subject to the relations
$$w(f)^*=w(-f),\quad w(f)w(g)=e^{-\2\kappa(f,g)}w(f+g).$$
The \emph{Fock state} $\omega_{P_1}$ over \CKK induced by $P_1$ is
again determined by condition \eqref{FS}. The GNS representation
$(\FF_{P_1},\pi_{P_1},\Omega_{P_1})$ of $\omega_{P_1}$ can be
identified with the representation $\pi_{P_1}$ given by formula
\eqref{PP1}, where $a^*(f)$ and $a(f)$, $f\in\KK_1$, now are Bosonic
creation and annihilation operators, acting on the symmetric Fock
space $\FF_{P_1}$ over $\KK_1$ with Fock vacuum vector
$\Omega_{P_1}$. All operators $\pi_{P_1}(a)$, $a\in\CKK$, are defined
on the invariant dense domain $\DD\subset\FF_{P_1}$ of algebraic
tensors, are closable, and fulfill
$\pi_{P_1}(a^*)\subset\pi_{P_1}(a)^*$.

The irreducible Fock representation of the Weyl algebra $\WK$ induced
by $P_1$ is obtained by identifying the Weyl operator $w(f)$,
$f=f^*\in\KK$, with the exponential of the closure of
$i\pi_{P_1}(f)$. The vacuum expectation value of $w(f)$ is
$$\langle{\Omega_{P_1}, w(f)\Omega_{P_1}}
\rangle=e^{-\frac{1}{4}\NORM{f}^2},\quad f=f^*.$$
Let \HH be a subspace of \KK with $\HH=\HH^*$, and let $\WW(\HH)$ be
the $C^*$-algebra generated by all $w(f)$ with $f=f^*\in\HH$. If
$\HH\SC$ is the orthogonal complement of \HH with respect to $\kappa$,
then
\begin{equation*}
  \WW(\HH)'=\WW(\HH\SC)''\qquad\text{(\emph{duality})}.
\end{equation*}

Every operator $V$ on \KK which preserves the form $\kappa$ and which
commutes with complex conjugation extends to a unique \emph{quasi-free
endomorphism} $\varrho_V$ of \CKK:
$$\varrho_V(f)=V(f),\quad f\in\KK,$$
and to a  unique *-endomorphism, denoted by the same symbol, of $\WK$:
$$\varrho_V(w(f))=w(V(f)),\quad f=f^*.$$
As shown in \cite{CB3}, an endomorphism $\varrho_V$ of $\WK$ can be
implemented by a Hilbert space of isometries $H_{\varrho_V}$ on
$\FF_{P_1}$ as in \eqref{IMP} if and only if the Hilbert--Schmidt
condition \eqref{IMPCOND} holds. Such operators $V$ form a semigroup
$$\END_{P_1}(\KK,\kappa)\DEF\SET{V\in\BB(\KK)}{V\+V=\1,\ \4{V}=V,\ [P_1,V]
\text{ is HS}}$$
isomorphic to the semigroup of all implementable quasi-free
endomorphisms  of $\WK$. Here we use the notation
$$A\+\DEF CA^*C$$
for the adjoint of an operator $A$ on \KK relative to the form
$\kappa$. Every $V\in\END_{P_1}(\KK,\kappa)$ has a well-defined
Fredholm index:
$$\IND V=\dim\ker V\+\in\{0,2,4,\dots,\infty\}.$$
The \emph{statistics dimension} $d_{\varrho_V}$ of $\varrho_V$ is in
the Bosonic case given by
\begin{equation}
	\label{ccr:DRV}
	d_{\varrho_V}\DEF\dim H_{\varrho_V}=
	\begin{cases}
  	1, & \IND V=0,\\
	  \infty, & \IND V\neq 0.
	\end{cases}
\end{equation}

In order to obtain an orthonormal basis in the Hilbert space
$H_{\varrho_V}\Omega_{P_1}$ associated with a fixed
$V\in\END_{P_1}(\KK)$, it is again convenient to extend the partial
Fock state $\omega_{P_1}\0\varrho_V^{-1}|_{\RAN\varrho_V}$ to a proper
Fock state $\omega_P$ which is unitarily equivalent to
$\omega_{P_1}$. The basis projection\footnote{In the Bosonic case, a
basis projection $P$ is an operator on \KK such that
$P=P^2=P\+=\1-\4{P}$, and such that $CP$ is positive definite on
$\RAN P$.} $P$ corresponding to this new Fock state has the form
$$P\DEF VP_1V\++p$$
where the basis projection $p$ of $\ker V\+$ is defined as follows
\cite{CB3}. Let $E$ be the orthogonal projection onto $\ker V\+$, and
let $A\DEF ECE$ be the operator describing the restriction of $\kappa$
to $\ker V\+$ with respect to the scalar product. Let $A_+$ be the
positive part of $A$, i.\,e.\ the unique positive operator such that
$A=A_+-\4{A_+}$ and $A_+\4{A_+}=0$. Let $A_+^{-1}$ be defined on $\RAN
A_+$ as the inverse of $A_+$, and on $(\RAN A_+)^\bot$ as zero. Then
$p$ is defined as
$$p\DEF A_+^{-1}C.$$

The class of all Fock states over $\WK$ which are unitarily equivalent
to $\omega_{P_1}$ is parameterized by symmetric\footnote{I.\,e.\
$T=\4{T^*}$.} Hilbert--Schmidt operators $T$ from $\KK_1$ to $\KK_2$
with $\NORM{T}<1$. The operator $T$ corresponding to $P$ is again
given by \eqref{TDEF}:
\begin{equation}
	\label{ccr:TDEF}
	T\DEF P_{21}{P_{11}}^{-1},
\end{equation}
whereas $P$ can be recovered from $T$ by
$$P=(P_1+T)(P_1+T\+T)^{-1}(P_1+T\+).$$
The cyclic vector $\Omega_P$ in $\FF_{P_1}$, unique up to a phase,
which induces the state $\omega_P$, is given by 
$$\Omega_P=\bigl({\det}(P_1+T\+T)\bigr)^{1/4}\exp(-\3\4{T}a^*a^*)
\Omega_{P_1}.$$ 
It belongs to $H_{\varrho_V}\Omega_{P_1}$. A complete orthonormal
basis in $H_{\varrho_V}\Omega_{P_1}$ is obtained by applying certain
isometries from the commutant of $\RAN\varrho_V$ to $\Omega_P$. The
basis projection $P$ leaves $\ker V\+$ invariant, and $\kappa$ is
positive definite on $\RAN P$. Let
\begin{equation}
	\label{ccr:K}
	\mathfrak{k}\DEF P(\ker V\+),\quad\text{with }
	\dim\mathfrak{k}=\3\IND V,
\end{equation}
and let $g_1,g_2,\dotsc$ be a basis in
$\mathfrak{k}$ such that $\kappa(g_j,g_k)=\delta_{jk}$. Let $\psi_j$
be the isometry obtained by polar decomposition of the closure of
$\pi_{P_1}(g_j)$. For any multi-index
$\alpha=(\alpha_1,\dots,\alpha_l)$ with
$1\leq\alpha_j\leq\alpha_{j+1}\leq\3\IND V$ ($\alpha=0$ if $l=0$), set
\begin{equation}
	\label{ccr:OA}
	\Omega_\alpha\DEF\psi_{\alpha_1}\dotsm\psi_{\alpha_l}\Omega_P.
\end{equation}
Then the $\Omega_\alpha$ form an orthonormal basis in
$H_{\varrho_V}\Omega_{P_1}$ and, by \eqref{HRV}, induce an orthonormal
basis in $H_{\varrho_V}$. The spaces $H_{\varrho_V}$ and
$H_{\varrho_V}\Omega_{P_1}$ can both be identified with the symmetric
Fock space over $\mathfrak{k}$ \cite{CB3}.
%%%%%%%%%%%%%%%%%%%%%%%%%%%%%%%%%%%%%%%%%%%%%%%%%%%%%%%%%%%%%%%%%
\subsection{Gauge invariant quasi-free endomorphisms (CCR)}
\label{sec:GICCR}
We assume again that a group $G\subset\END_{P_1}(\KK,\kappa)$
consisting of unitary operators $U$ which commute with $P_1$ (so that
$G$ can be identified with a subgroup of $U(\KK_1)$) acts by
second quantization $\Gamma(U)$ on $\FF_{P_1}$. \emph{Gauge invariant}
implementable quasi-free endomorphisms correspond to the elements of
the semigroup
$$\END_{P_1}(\KK,\kappa)^G\DEF\SET{V\in\END_{P_1}
(\KK,\kappa)}{[V,G]=0}.$$
To determine the representation $D_{\varrho_V}$ of $G$ on the Hilbert
space $H_{\varrho_V}$ associated with $V\in\END_{P_1}(\KK,\kappa)^G$,
it suffices to consider the action of $\Gamma(U)$ on the vectors
$\Omega_\alpha$ defined in \eqref{ccr:OA}. In contrast to the
Fermionic case (cf.\ Lemma~\ref{lem:GP}), there is no simple
transformation law for the $\psi_j$ under $G$. They obey however a
linear transformation law when applied to $\Omega_P$; in fact, one can
show that $\Omega_\alpha$ is proportional to
\begin{equation}
	\label{ccr:OA'}
	\pi_{P_1}(g_{\alpha_1})\dotsm\pi_{P_1}(g_{\alpha_l})
	\exp(-\3\4{T}a^*a^*)\Omega_{P_1}
\end{equation}
(cf.\ \cite{CB4}; taking the closures of the $\pi_{P_1}(g_j)$ is
tacitly assumed here). The behavior of the $\pi_{P_1}(g_j)$ under
gauge transformations is obvious.
\begin{Lem}
  \label{ccr:lem:GEXP}
  Let $V\in\END_{P_1}(\KK,\kappa)^G$ be given, and let $T$ be defined
  by \eqref{ccr:TDEF}. Then $\exp(-\3\4{T}a^*a^*)\Omega_{P_1}$ is
  invariant under all gauge transformations $\Gamma(U)$, $U\in G$.
\end{Lem}
\begin{proof}
  Let $E$ be the orthogonal projection onto $\ker V\+=C\ker V^*$, and
  let $A\DEF ECE$ as in Section~\ref{sec:PRECCR}. Then $E$ and $A$
  commute with $G$ because $V$ and $P_1$ do so. Therefore the positive
  part $A_+$ of $A$ and the operator ${A_+}^{-1}$ defined in
  Section~\ref{sec:PRECCR} also commute with $G$. It follows that
  $P=VP_1V\++{A_+}^{-1}C$ and $T\DEF P_{21}{P_{11}}^{-1}$ commute with
  $G$ as well. 
  
  Arguing as in the proof of Lemma~\ref{lem:EXP}, one finds for $U\in
  G$
  $$\Gamma(U)\bigl(-\3\4{T}a^*a^*\bigr)^n\Omega_{P_1}
  =\bigl(-\3(U\4{T}U\+)a^*a^*\bigr)^n\Omega_{P_1}$$ 
  and finally
  $$\Gamma(U)\exp\bigl(-\3\4{T}a^*a^*\bigr)\Omega_{P_1}=
  \exp\bigl(-\3(U\4{T}U\+)a^*a^*\bigr)\Omega_{P_1}=
  \exp\bigl(-\3\4{T}a^*a^*\bigr)\Omega_{P_1}.$$
\end{proof}
\begin{Thm}
  \label{ccr:thm:CHARGE}
  Let $P_1$ be a basis projection of $(\KK,\kappa)$, let $G$ be a
  group of unitary operators on \KK commuting with $P_1$ and with
  complex conjugation, and let $V\in\END_{P_1}(\KK,\kappa)^G$.  Then
  the subspace $\mathfrak{k}$ defined in \eqref{ccr:K} is invariant
  under $G$, and the unitary representation $D_{\varrho_V}$ of $G$ on
  the Hilbert space of isometries $H_{\varrho_V}$ which implements
  $\varrho_V$ in the Fock representation determined by $P_1$ is
  unitarily equivalent to the representation
  $\mathfrak{S}_\mathfrak{k}$ on the symmetric Fock space over
  $\mathfrak{k}$ that is obtained by taking symmetric 
  tensor powers of the representation on $\mathfrak{k}$:
  \begin{equation} 
		\label{ccr:REPG}
  	D_{\varrho_V}\simeq\mathfrak{S}_{\mathfrak{k}}.
  \end{equation}
\end{Thm}
\begin{proof}
  $\mathfrak{k}$ is invariant under $G$ because $\ker V\+$ is
  invariant and because $P$ commutes with $G$ (see the proof of
  Lemma~\ref{ccr:lem:GEXP}). The assertion hence follows from
  \eqref{ccr:OA'} and Lemma~\ref{ccr:lem:GEXP}.
\end{proof}
\begin{Rms}
	(i) Theorem~\ref{ccr:thm:CHARGE} shows that non-surjective quasi-free
	endomorphisms of the CCR algebra are even ``more reducible'' than
	endomorphisms of the CAR algebra in that they are always
	\emph{infinite} direct sums, a fact which explains the generic
	occurrence of infinite statistics in the CCR case (cf.\
	\eqref{ccr:DRV}). Again, each closed subspace of
	$H_{\varrho_V}\Omega_{P_1}$ spanned by the $\Omega_\alpha$ with
	length of $\alpha$ fixed is invariant under $G$.

	(ii) Any representation of $G$ which is contained in $\KK_1$ with
	infinite multiplicity is realized on some space $\mathfrak{k}$
	belonging to a $V\in\END_{P_1}(\KK,\kappa)^G$.

	(iii) Quasi-free automorphisms are less interesting in the CCR
	case because they are all neutral: $D_{\varrho_V}$ is the trivial
	representation of $G$ if $\IND V=0$.

	(iv) If the single-particle space $\KK_1$ has a
	particle-antiparticle symmetry, then every
	$V\in\END_{P_1}(\KK,\kappa)^G$ has a conjugate in
	$\END_{P_1}(\KK,\kappa)^G$, just as in the Fermionic case. 
\end{Rms}
%%%%%%%%%%%%%%%%%%%%%%%%%%%%%%%%%%%%%%%%%%%%%%%%%%%%%%%%%%%%%%%
\section{Concluding Remarks}
As we have seen, gauge invariant implementable quasi-free
endomorphisms of the CAR and CCR algebras with statistics dimension
$d\neq1$ restrict to reducible endomorphisms of the observable
algebra.  In typical cases, e.g.\ if $G$ is isomorphic to one of the
classical compact Lie groups, any irreducible representation of the
group is equivalent to a subrepresentation of some tensor power of the
defining representation.  In such cases there will exist quasi-free
endomorphisms, behaving like ``master endomorphisms'', which contain
each superselection sector as a subrepresentation. 

It is an interesting problem how to obtain the irreducible
``subobjects'' of a quasi-free endomorphism $\varrho$. Suppose that
$\{\Psi_j\}$ is an (incomplete) orthonormal set in $H_\varrho$ which
transforms irreducibly under $G$. According to the general theory
\cite{DHR3}, there should exist a gauge invariant isometry $\Phi$ on
Fock space with $\RAN\Phi=\oplus_j\RAN\Psi_j$.  The corresponding
irreducible endomorphism $\varrho_\Phi$ (which is not quasi-free)
would then be given by
$$\varrho_\Phi(a)\DEF\Phi^*\Bigl(\sum_j\Psi_ja\Psi_j^*\Bigr)\Phi.$$
Collections of gauge invariant isometries $\{\Phi_j\}$ fulfilling the
Cuntz relations would permit to define direct sums of quasi-free
endomorphisms $\{\varrho_j\}$:
\begin{equation*}
  (\oplus_j\varrho_j)(a)\DEF\sum_j\Phi_j\varrho_j(a)\Phi_j^*,
\end{equation*}
so that one would get the whole Doplicher--Roberts category generated
by quasi-free endomorphisms.

Another important question is how to find basis-independent examples
of, say, localized isometries $v$ with index one on the
single-particle space \HH (see \eqref{HL2}) of the time-zero Dirac
field, such that the implementability condition \eqref{Ev}
holds. Recall that our construction of such an operator in
Eq.~\eqref{v} made essential use of the existence of local Fourier
bases on the circle. Of particular interest would be the massive case
in two dimensions, where one might hope to find localized quasi-free
endomorphisms obeying non-Abelian braid group statistics. However,
preliminary calculations based on the explicit formulas in
\cite{CB1,CB3} indicate that the commutation relations of implementers
corresponding to irreducible subobjects of quasi-free endomorphisms
only admit Abelian braid group statistics.
\begin{acknowledgements}
	The author profited from discussions with P.~Grinevich on the
  estimates mentioned in the proof of Lemma~\ref{lem:HS}.
\end{acknowledgements}
%%%%%%%%%%%%%%%%%%%%%%%%%%%%%%%%%%%%%%%%%%%%%%%%%%%%%%%%%%%%%%%%%
\providecommand{\bysame}{\leavevmode\hbox to3em{\hrulefill}\thinspace}

\end{document}